\documentclass[11pt]{amsart}
\usepackage{amssymb}

\numberwithin{equation}{section}

\theoremstyle{plain} 
\newtheorem{theor}[equation]{Theorem}
\newtheorem{cor}[equation]{Corollary}
\newtheorem{lem}[equation]{Lemma}
\newtheorem{proposition}[equation]{Proposition}
\newtheorem{resultat}[equation]{}

\theoremstyle{definition}
\newtheorem{defin}[equation]{Definition}
\newtheorem{numerotation}[equation]{}

\theoremstyle{remark}

\newtheorem{ex}[equation]{Example}

\newcommand{\Character}{\chi}
\newcommand{\lineaire}{\Theta}

\def\build#1_#2^#3{\mathrel{\mathop{\kern0pt#1}\limits_{#2}^{#3}}}

\begin{document}
\title{\bf Batalin-Vilkovisky algebras and cyclic cohomology of Hopf algebras}
\author{Luc Menichi}
\address{UMR 6093 associ\'ee au CNRS\\
Universit\'e d'Angers, Facult\'e des Sciences\\
2 Boulevard Lavoisier\\49045 Angers, FRANCE}
\email{Luc.Menichi@univ-angers.fr}

\begin{abstract}
We show that the Connes-Moscovici cyclic cohomology of a Hopf algebra
equipped with a character has a Lie bracket of degree $-2$.
More generally, we show that a "cyclic operad with multiplication"
is a cocyclic module whose cohomology is a Batalin-Vilkovisky algebra
and whose cyclic cohomology is a graded Lie algebra of degree $-2$.
This explain why the Hochschild cohomology algebra of a symmetric
algebra is a Batalin-Vilkovisky algebra.
\end{abstract}
\maketitle

\section{Introduction}
Let $\Bbbk$ be an arbitrary commutative ring
and denote by $\mathcal{H}$ an (ungraded) bialgebra over $\Bbbk$.
We denote by $\Omega\mathcal{H}$ the Adams Cobar construction
on $\mathcal{H}$. Its cohomology is $\text{Cotor}^{*}_\mathcal{H}(\Bbbk,\Bbbk)$.
It results from~\cite[p. 65]{Gerstenhaber-Schack:algbqgad}
that $\text{Cotor}^{*}_\mathcal{H}(\Bbbk,\Bbbk)$
has a Gerstenhaber algebra structure.

On the other hand, assume that $\mathcal{H}$ has an involutive
antipode or more generally that $\mathcal{H}$ is an Hopf
algebra equipped with a modular pair in involution
where the group like element is the unit $1$ of $\mathcal{H}$.
Connes and Moscovici~\cite{Connes-Mosco:hopfacctit,Connes-Mosco:cyclchas}
have proved that $\Omega\mathcal{H}$ has a canonical cocyclic module
structure.
Since a cocyclic structure defines a Connes coboundary map $B$ in
cohomology and since a Batalin-Vilkovisky algebra
(Definition~\ref{definition BV algebre})
 is a Gerstenhaber algebra equipped
with an operator $B$,
it is natural to conjecture that 
$\text{Cotor}^{*}_\mathcal{H}(\Bbbk,\Bbbk)$ is
a Batalin-Vilkovisky algebra. The first result of this paper
is to prove that conjecture.
\begin{theor}\label{Cobar algebre de Hopf}
Let $\mathcal{H}$ be an Hopf algebra endowed with a modular pair in
involution $(\Character,1)$. Then

a) The canonical algebra structure of the Cobar construction
on $\mathcal{H}$ together with its Connes-Moscovici cocyclic structure,
defines a Batalin-Vilkovisky
algebra structure on $\text{Cotor}^{*}_{\mathcal{H}}(\Bbbk,\Bbbk)$.

b) The cyclic cohomology of $\mathcal{H}$, denoted
$HC^*_{(\Character,1)}(\mathcal{H})$, is a graded Lie algebra of degree $-2$.
\end{theor}
The easiest way to see that the cotorsion product of a bialgebra
$\mathcal{H}$ is a Gerstenhaber algebra, is to remark as
in~\cite{Gerstenhaber-Schack:algbqgad} that the Cobar construction on $\mathcal{H}$,
$\Omega\mathcal{H}$, is an operad with multiplication
 (Definition~\ref{definition operade avec multiplication})
and to apply the following general theorem.
\begin{resultat}\cite{Gerstenhaber-Schack:algbqgad,Gerstenhaber-Voronov:hgamso,McClure-Smith:deligneconj}\label{Theoreme Gerstenhaber}
\noindent a) Each operad with multiplication $O$ is a cosimplicial
module.
Denote by $\mathcal{C}^{*}(O)$ the associated cochain complex. 

\noindent b) Its cohomology $H^*(\mathcal{C}^{*}(O))$
is a Gerstenhaber algebra.
\end{resultat}
To prove Theorem~\ref{Cobar algebre de Hopf}, we proceed similarly:

-we introduce the notion of cyclic operad with multiplication
(Definition~\ref{operade cyclique avec multiplication}),

-we show in section~\ref{cyclic cohomology of Hopf algebras}
that $\Omega\mathcal{H}$ is a cyclic operad with multiplication.

-we prove the main result of this paper. 
\begin{theor}\label{Theoreme BV algebre}
If $O$ is a cyclic operad with a multiplication then

a) the structure of cosimplicial module on $O$ extends to a structure
of cocyclic module,

b) the Connes coboundary map $B$ on $\mathcal{C}^{*}(O)$
induces a natural structure of Batalin-Vilkovisky algebra
on the Gerstenhaber algebra $H^*(\mathcal{C}^{*}(O))$ and

c) the cyclic cohomology of $\mathcal{C}^{*}(O)$, $HC^*(\mathcal{C}^{*}(O))$,
has naturally a graded Lie algebra structure of degree $-2$.
\end{theor}
Part b) of Theorem~\ref{Theoreme BV algebre} is inspired
by a result (See section~\ref{resultat de McClure and Smith})
announced by McClure and Smith~\cite{McClure-Smith:slidesnorthwestern}.

As a second application of Theorem~\ref{Theoreme BV algebre},
we show
\begin{theor}\label{Hochschild algebre symmetrique}
Let $A$ be an algebra equipped with an isomorphism of
$A$-bimodules $\lineaire:A\buildrel{\cong}\over\rightarrow A^{\vee}$
(i. e. $A$ is a symmetric algebra).
Then

a) the Connes coboundary map on $HH^{*}(A,A^{\vee})$ defines via
the isomorphism $HH^{*}(A,\lineaire):HH^{*}(A,A)
\buildrel{\cong}\over\rightarrow HH^{*}(A,A^{\vee})$
a structure of Batalin-Vilkovisky algebra on the Gerstenhaber algebra
$HH^{*}(A,A)$.

b) The cyclic cohomology of $A$, $HC^*(A)$ is a
graded Lie algebra of degree $-2$.
\end{theor}

Part a) of Theorem~\ref{Hochschild algebre symmetrique}
has been proved by Tradler~\cite{Tradler:bvalgcohiiip}.
In fact, he proved part a) much more generally,
for "homotopy" symmetric algebras.
Our proof for "strict" symmetric algebras is much simpler.

Tamarkin and Tsygan~\cite[Conjecture 0.13]{Tamarkin-Tsygan:ncdchBValg}
have conjectured a related result at the chain level.
See also McClure and Smith~\cite{McClure-Smith:slidesnorthwestern}.
Moreover Tamarkin and Tsygan have mentionned a relation
between Part a) of Theorem~\ref{Hochschild algebre symmetrique}
and Connes-Moscovici cyclic cohomology of Hopf
algebras~\cite{Tamarkin-Tsygan:ncdchBValg}.
Theorem~\ref{Theoreme BV algebre} establishes such relation.

The main tools in the proof of Theorem~\ref{Theoreme BV algebre}
are the following results which have their own interest.
Let $\overline{\mathcal{C}}^*(O)$ be the normalized cochain complex
associated to the cyclic operad with multiplication $O$
and $B$ the Connes normalized coboundary map on
$\overline{\mathcal{C}}^*(O)$.
Denote by $\cup$ the cup product and by $\overline{\circ}$
the composition product in $\overline{\mathcal{C}}^*(O)$
(See (\ref{definition du cup}) and (\ref{definition du composition})).
\begin{lem}\label{multiple de deux}
There is a bilinear map Z (See (\ref{definition de Z}))
of degree $-1$ such that
$$B(f\cup g)=Z(f,g)+(-1)^{mn}Z(g,f),\quad
\forall f\in\overline{\mathcal{C}}^{m}(O),
g\in\overline{\mathcal{C}}^{n}(O).$$
\end{lem}
\begin{proposition}\label{moitie BV}
There is a bilinear map $H$ (See (\ref{definition de H}))
of degree $-2$ such that, for any $f\in\overline{\mathcal{C}}^{m}(O)$
and $g\in\overline{\mathcal{C}}^{n}(O)$,
\begin{multline*}
(-1)^{m}\left(Z(f,g)-(Bf)\cup g\right)-f\overline{\circ}g\\
=dH(f,g)+H(df,g)+(-1)^{m-1}H(f,dg).
\end{multline*}
\end{proposition}
We give now the plan of the paper:

{\bf 2) operads with multiplication.}
This section is a review on operads with multiplication.
We recall the definition of operad with multiplication. We define
the structure of Gerstenhaber algebra associated to an operad with
multiplication. We recall the two fundamental examples of operad with
multiplication:

-the endomorphism operad of an algebra,

-the Cobar construction on a bialgebra.

{\bf 3) cyclic operad with multiplication.}
We introduce the notions of cyclic (non-$\Sigma$) operad
and cyclic operad with multiplication.

{\bf 4) Hochschild cohomology of an symmetric algebra.}
We prove Theorem~\ref{Hochschild algebre symmetrique}
by showing that the endomorphism operad of a symmetric
algebra is a cyclic operad with multiplication.

{\bf 5) Cyclic cohomology of Hopf algebras.}
We recall what is an Hopf algebra $\mathcal{H}$ endowed with a modular pair in
involution of the form $(\Character,1)$ and we prove
Theorem~\ref{Cobar algebre de Hopf} by showing that the Cobar construction
on $\mathcal{H}$ is a cyclic operad with multiplication.

{\bf 6) Proof of parts a) and b) of Theorem~\ref{Theoreme BV algebre}.}
We prove part a) of Theorem~\ref{Theoreme BV algebre} and
Lemma~\ref{multiple de deux}. Then we deduce part b) from
Lemma~\ref{multiple de deux} and Proposition~\ref{moitie BV}.
Finally, we prove Proposition~\ref{moitie BV}.

{\bf 7) Proof of part c) of Theorem~\ref{Theoreme BV algebre}.}
We define the Lie bracket on cyclic cohomology
in the same way as Chas and Sullivan define a Lie bracket
on $S^{1}$-equivariant homology in~\cite{Chas-Sullivan:stringtop}.

{\bf 8) Comparison with McClure and Smith.}
We compare Theorem~\ref{Theoreme BV algebre} with two results
announced by McClure and Smith~\cite{McClure-Smith:slidesnorthwestern}.

{\em Acknowledgment:} We wish to thank Jean-Claude Thomas for his
constant support.
\section{operads with multiplication}
\begin{numerotation}
A {\it Gerstenhaber algebra} is a
graded module $G=\{G^i\}_{i\in
\mathbb{Z}}$
equipped with two linear maps
$$
\cup:G^i \otimes G^j \to G^{i+j} \,, \quad x\otimes y \mapsto x\cup y
$$
$$
\{-,-\}:G^i \otimes G^j \to G^{i+j-1} \,, \quad x\otimes y \mapsto \{x,y\}
$$

such that:

\noindent a) the cup product $\cup$ makes G into a graded commutative
algebra

\noindent b) the bracket $\{-,-\}$ gives $G$ a structure of graded
Lie algebra of degree $-1$. This means that for each $a$, $b$ and $c\in G$

$\{a,b\}=-(-1)^{(\vert a\vert-1)(\vert b\vert-1)}\{b,a\}$ and 

$\{a,\{b,c\}\}=\{\{a,b\},c\}+(-1)^{(\vert a\vert-1)(\vert b\vert-1)}
\{b,\{a,c\}\}.$

\noindent c)  the cup product and the Lie bracket satisfy the Poisson rule.
This means that for any  $a  \in G^k$  
the adjunction map $\{a,-\}:G^i \to G^{i+k-1} \,, \quad
b\mapsto \{a,b\}$ is  a
$(k+1)$-derivation: ie. for $a$, $b$, $c\in G$,
$\{ a, b c\} = \{ a,b\} c + (-1)^{\vert b\vert (\vert a\vert -1)} 
b \{ a,c\}$.
\end{numerotation}

Usually, this definition is given for a lower graded module $G=\{G_i\}_{i\in\mathbb{Z}}$. If you put $G_i=G^{-i}$ as usual, you pass from an upper degree graded module to a lower graded module and the Lie bracket is of the usual (lower)
degree $+1$.

\begin{numerotation}
In this paper, {\it operad} means (non-$\Sigma$) operad in the category of
$\Bbbk$-modules. That is:
a sequence of modules $\{O(n)\}_{n\in\mathbb{N}}$, an identity element $id\in O(1)$
and structure maps

$\gamma:O(n)\otimes O(i_1)\otimes\dots\otimes O(i_n)\rightarrow O(i_1+\dots+i_n)$

$f\otimes g_1\otimes\dots\otimes g_n\mapsto \gamma(f;g_1,\dots,g_n)$

satisfying associativity and unit~\cite{Markl-Shnider-Stasheff:opeatp}.
\end{numerotation}
Hereaffter we use mainly the composition operations
$\circ_i:O(m)\otimes O(n)\rightarrow O(m+n-1)$
$f\otimes g\mapsto f\circ_i g$
defined for $m\in\mathbb{N}^{*}$, $n\in\mathbb{N}$ and $1\leq i\leq m$
by $f\circ_i g:=\gamma(f;id,\dots,g,id,\dots,id)$ where $g$ is the $i$-th element
after the semicolon.

\begin{ex}~\cite{Markl-Shnider-Stasheff:opeatp}
Let $V$ be a module. The {\it endomorphism operad} of $V$ is the operad
$\mathcal{E}nd_V$ defined by $\mathcal{E}nd_V(n):=\text{Hom}(V^{\otimes n},V)$.
The identity element of $\mathcal{E}nd_V$ is the identity map $id_V:V\rightarrow V$.
\end{ex}
\begin{numerotation}\label{definition operade avec multiplication}
An {\it operad with multiplication} is a operad equipped with an element $\mu\in O(2)$ called the multiplication and an element $e\in O(0)$
such that
$\mu\circ_1\mu=\mu\circ_2\mu$
and $\mu\circ_1 e=id=\mu\circ_2 e.$
\end{numerotation}
In~\cite{Gerstenhaber-Schack:algbqgad}, an operad with multiplication
is called a strict unital comp algebra.

Let $Ass$ be the (non-$\Sigma$) associative operad: $Ass(n):={\Bbbk}$.
An operad $O$ is an operad with multiplication if and only if
$O$ is equipped with a morphism of operads $Ass\rightarrow O$.
\begin{proof}[Sketch of proof of~\ref{Theoreme Gerstenhaber}]

a) The coface maps $\delta_{i}:O(n)\rightarrow O(n+1)$ and codegeneracy maps
$\sigma_{j}:O(n)\rightarrow O(n-1)$
are defined \cite{McClure-Smith:deligneconj} by
$\delta_{0}f=\mu\circ_2 f$,
$\delta_{i}f=f\circ_i\mu$,
$\delta_{n+1}f=\mu\circ_1 f$ and $\sigma_{i-1}f=f\circ_{i}e$ for $1\leq i\leq n$.

b) The associated cochain complex $\mathcal{C}^{*}(O)$ is the cochain
complex whose differential $d$ is given by $$d:=\sum_{i=0}^{n+1}(-1)^{i}\delta_{i}:O(n)\rightarrow O(n+1).$$
The linear maps $\cup:O(m)\otimes O(n)\rightarrow O(m+n)$ defined by
\begin{equation}\label{definition du cup}
f\cup g:=(\mu\circ_1 f)\circ_{m+1}g=(\mu\circ_2 g)\circ_{1}f
\end{equation}
gives $\mathcal{C}^{*}(O)$ a structure of differential graded algebra.
The linear maps of degree $-1$
$$\overline{\circ},\{-,-\}:O(m)\otimes O(n)\rightarrow O(m+n-1)$$
are defined by
\begin{equation}\label{definition du composition}
f\overline{\circ}g:=(-1)^{(m-1)(n-1)}\sum_{i=1}^{m}
(-1)^{(n-1)(i-1)}f\circ_i g
\end{equation}
and $$\{f,g\}:=f\overline{\circ}g-(-1)^{(m-1)(n-1)}g\overline{\circ}f.$$
The bracket $\{-,-\}$ defines a structure of differential graded Lie algebra
of degree $-1$ on $\mathcal{C}^{*}(O)$.
After passing to cohomology, the cup product $\cup$ and the bracket $\{-,-\}$
satisfy the Poisson rule.
\end{proof}
\begin{cor}\cite{Gerstenhaber:cohosar}\label{operad d'endomorphisme d'une algebre}
The Hochschild cohomology of an algebra is a Gerstenhaber algebra.
\end{cor}
\begin{proof}
Let $A$ be an associative algebra with multiplication
$\mu:A\otimes A\rightarrow A$ and unit $e:\Bbbk\rightarrow A$.
Then the endomorphism operad $\mathcal{E}nd_A$ of $A$ equipped
with $\mu$ and $e$ is an operad with multiplication.
The Hochschild cochain complex of $A$, denoted $\mathcal{C}^{*}(A,A)$,
is the cochain complex $\mathcal{C}^{*}(\mathcal{E}nd_A)$ associated
to the endomorphism operad of $A$.
\end{proof}
\begin{cor}\cite[p. 65]{Gerstenhaber-Schack:algbqgad}\label{Cobar d'une bialgebre}
Let $\mathcal{H}$ be a bialgebra. Then
$\text{Cotor}^{*}_\mathcal{H}(\Bbbk,\Bbbk)$ is a Gerstenhaber algebra.
\end{cor}
\begin{proof}
Denote by $\mu$ and $1$ the multiplication and the unit of $\mathcal{H}$.
Denote by $\Delta$ and $\varepsilon$ the diagonal and the counit of
$\mathcal{H}$.
For each $n\in\mathbb{N}$, denote by
$\Delta^{n-1}: \mathcal{H}\rightarrow \mathcal{H}^{\otimes n}$
the $(n-1)$ iterated diagonal defined by
$\Delta^{-1}:=\varepsilon$, $\Delta^{0}:=Id_\mathcal{H}$ and
$\Delta^{n+1}:=(\Delta\otimes id_\mathcal{H}^{\otimes n})\circ\Delta^{n}$.
For an element $a\in\mathcal{H}$, we denote
$\Delta^{n-1}a:=a^{(1)}\otimes\dots\otimes a^{(n)}$ or simply
$a^{1}\otimes\dots\otimes a^{n}$. Here the sum is implicit and contrarly to
Sweedler notation, we use upperscripts instead of lowerscripts, since we will
need indices but no powers.

Consider the operad with multiplication $O$ defined by
$O(n):=\mathcal{H}^{\otimes n}$ and if
$a_1\otimes\dots\otimes a_m\in\mathcal{H}^{\otimes m}$ and
$b_1\otimes\dots\otimes b_n\in\mathcal{H}^{\otimes n}$,

$$(a_1\otimes\dots\otimes a_m)\circ_i(b_1\otimes\dots\otimes b_n):=$$
$$
a_1\otimes\dots\otimes a_{i-1}\otimes(\Delta^{n-1}a_i)\cdot(b_1\otimes\dots\otimes b_n)\otimes a_{i+1}\otimes\dots\otimes a_m=$$
$$a_1\otimes\dots\otimes a_{i-1}\otimes a_i^{1}b_1\otimes\dots\otimes a_i^{n}b_n\otimes a_{i+1}\otimes\dots\otimes a_m.$$
Here $\cdot$ is the product obtained by tensorization in $\mathcal{H}^{\otimes n}$.
The identity element $id$ of $O$ is $1\in\mathcal{H}^{\otimes 1}$.
The multiplication $\mu$ is $1\otimes 1\in\mathcal{H}^{\otimes 2}$.
The element $e$ of $O$ is the unit of $\Bbbk$,
$1_\Bbbk\in\mathcal{H}^{\otimes 0}$.
The cochain complex associated to this operad is the Cobar construction
on $\mathcal{H}$, denoted usually $\Omega\mathcal{H}$.
Since
$\text{Cotor}^*_\mathcal{H}(\Bbbk,\Bbbk)=H^*(\Omega\mathcal{H})$,
the result follows from~\ref{Theoreme Gerstenhaber}.
Remark that if $\mathcal{H}$ is not cocommutative,
$O$ is not in general a symmetric operad.
\end{proof}
\section{Cyclic operads with multiplication}
\begin{numerotation}\label{definition BV algebre}
A {\it Batalin-Vilkovisky algebra} is a Gerstenhaber algebra
$G$ equipped with a degree $-1$ linear map $B:G^{i}\rightarrow G^{i-1}$
such that $B\circ B=0$ and
\begin{equation}\label{relation BV} 
\{a,b\}=(-1)^{\vert a\vert}\left(B(a\cup b)-(B a)\cup b-(-1)^{\vert
  a\vert}a\cup(B b)\right)
\end{equation}
for $a$ and $b\in G$.
\end{numerotation}
\begin{defin}
A {\it cyclic (non-$\Sigma$) operad} is a (non-$\Sigma$) operad $O$
equipped with linear maps $\tau_n:O(n)\rightarrow O(n)$
for $n\in\mathbb{N}$ such that
\begin{equation}
\forall n\in\mathbb{N},\quad\tau_n^{n+1}=id_{O(n)},
\end{equation}
\begin{equation}\label{deux cyclic}
 \forall m\geq 1,n\geq 1,
\quad\tau_{m+n-1}(f\circ_1 g)= \tau_n g\circ_n \tau_m f,
\end{equation}
\begin{equation}\label{trois cyclic}
\forall m\geq 2, n\geq 0, 2\leq i\leq m,
\quad\tau_{m+n-1}(f\circ_i g)=\tau_m f\circ_{i-1}g,
\end{equation}
for each $f\in O(m)$ and $g\in O(n)$. In particular, we have
$\tau_1 id=id$.
\end{defin}
This definition is taken from~\cite[p. 247-8]{Markl-Shnider-Stasheff:opeatp}
except that since our operad $O$ is not necessarly symmetric, we don't
assume that the action of the cyclic group $\mathbb{Z}/(n+1)\mathbb{Z}$
on $O(n)$
extends to an action of the symmetric group of order $n+1$, $S_{n+1}$.

Remark that~(\ref{deux cyclic}) and~(\ref{trois cyclic}) are equivalent to
\begin{equation}\label{deux duale}
\forall m\geq 1,n\geq 1,
\quad\tau_{m+n-1}^{-1}(f\circ_m g)= \tau_n^{-1} g\circ_1 \tau_m^{-1} f,
\end{equation}
\begin{equation}\label{trois duale}
\forall m\geq 2, n\geq 0, 1\leq i\leq m-1,
\quad\tau_{m+n-1}^{-1}(f\circ_i g)=\tau_m^{-1} f\circ_{i+1}g,
\end{equation}

If instead of~(\ref{deux cyclic}) and~(\ref{trois cyclic}), $\tau_n$
satisfies
\begin{align*}
\tau_{m+n-1}(f\circ_m g)= \tau_n g\circ_1 \tau_m f,\\
\tau_{m+n-1}(f\circ_i g)=\tau_m f\circ_{i+1}g,
\end{align*}
like in the original definition of cyclic operad of Getzler and
Kapranov~\cite[(2.2)]{Getzler-Kapranov:cycloch},
replace $\tau_n$ by $\tau_n^{-1}$.

We will use the following generalizations of~(\ref{trois duale})
and of~(\ref{deux duale}):
For each $m\geq 1$, $n\geq 0$, $1\leq i\leq m$ and $j\in\mathbb{Z}$,

\begin{equation}\label{trois generalise}
\text{ if }1\leq i+j\leq m\text{ then }
\tau_{m+n-1}^{-j}(f\circ_i g)=
\tau_m^{-j} f\circ_{i+j}g,
\end{equation}
\begin{multline}\label{deux generalise}
\text{ if }m+1\leq i+j\leq m+n\text{ then }\\
\tau_{m+n-1}^{-j}(f\circ_i g)=
\tau_n^{-j+m-i} g\circ_{i+j-m}\tau_m^{i-m-1} f.
\end{multline}

\begin{defin}\label{operade cyclique avec multiplication}
A {\it cyclic operad with multiplication} is an operad
which is both an operad with multiplication and a cyclic operad
such that 
\begin{equation*}
\tau_2\mu=\mu.
\end{equation*}
\end{defin}
The operad $Ass$ is a cyclic operad: the cyclic group
$\mathbb{Z}/(n+1)\mathbb{Z}$
act trivially on $Ass(n):=\Bbbk$.
A cyclic operad $O$ is an cyclic operad with multiplication if and only if
$O$ is equipped with a morphism of cyclic operads $Ass\rightarrow O$.
\section{Hochschild cohomology of a symmetric algebra}
{\bf The cyclic endomorphism operad~\cite{Markl-Shnider-Stasheff:opeatp}.}
Let $V$ be a module equipped with a
bilinear form $\varphi:V\otimes V\rightarrow\Bbbk$ such that the associated
right linear map
$\lineaire:V\buildrel{\cong}\over\rightarrow V^{\vee}$,
$v\mapsto \varphi(-,v)$, is an isomorphism.
Consider the adjonction map
\begin{equation}\label{adjunction map}
Ad:\text{Hom}(V^{\otimes n},V^{\vee})\buildrel{\cong}\over\rightarrow
\text{Hom}(V^{\otimes n+1},\Bbbk)
\end{equation}
 which associates to any
$g\in\text{Hom}(V^{\otimes n},V^{\vee})$, the map
$$Ad(g):V^{\otimes n+1}\rightarrow\Bbbk,\quad
v_0,v_1,\cdots,v_n\mapsto g(v_1,\cdots,v_n)(v_0).$$
The composite
$$
\text{Hom}(V^{\otimes n},V)
\buildrel{\text{Hom}(V^{\otimes n},\lineaire)}\over\longrightarrow
\text{Hom}(V^{\otimes n},V^{\vee})
\buildrel{Ad}\over\rightarrow
\text{Hom}(V^{\otimes n+1},\Bbbk)
$$
is an isomorphism. Explicitly this composite send
$f\in\text{Hom}(V^{\otimes n},V)$ to the linear map
$\widehat{f}:V^{\otimes n+1}\rightarrow\Bbbk $ defined
by
$$\widehat{f}(v_0,v_1,\cdots,v_n)=\varphi(v_0,f(v_1,\cdots,v_n))\quad
\text{for }v_0, v_1, \cdots, v_n\in V.$$
The cyclic group $\mathbb{Z}/(n+1)\mathbb{Z}$ acts on $V^{\otimes n+1}$
by permutations of factors:
\begin{equation}\label{cyclic operator}
t_n(v_0,\cdots,v_n):=(v_n,v_0,\cdots,v_{n-1})\quad
\text{for }(v_0,\cdots,v_n)\in V^{\otimes n+1}.
\end{equation}
Define $\tau_n:=t_n^{\vee}:\text{Hom}(V^{\otimes n+1},\Bbbk)
\rightarrow \text{Hom}(V^{\otimes n+1},\Bbbk)$.
Using the identification $f\mapsto\widehat f$, we define
$\tau_n:\text{Hom}(V^{\otimes n},V)\rightarrow
\text{Hom}(V^{\otimes n},V)$ by
$\widehat{\tau_n f}:=\tau_n\widehat{f}$ for $f\in\text{Hom}(V^{\otimes
  n},V)$.
Explicitly, $\tau_n(f)$ is the unique map such that
$$
\varphi(v_0,\tau_n f(v_1,\cdots,v_n))=
\varphi(v_n,f(v_0,\cdots,v_{n-1}))\quad
\text{for } v_0,\cdots,v_n\in V.$$
The endomorphism operad of $V$, $\mathcal{E}nd_V$, equipped with this last
linear map $\tau_n:\mathcal{E}nd_V(n)\rightarrow\mathcal{E}nd_V(n)$
is a cyclic operad if and only if the bilinear form $\varphi$
is symmetric.

{\bf Hochschild (co)homology.}
Let $A$ be an algebra.
Let $M$ be an $A$-bimodule.
Denote by $\mathcal{C}^*(A,M)$ the Hochschild
cochain complex of $A$ with coefficient in $M$~\cite[1.5.1]{LodayJ.:cych}
and by $\mathcal{C}_*(A,M)$ the Hochschild chain
complex~\cite[1.1.1]{LodayJ.:cych}.
Recall that $\mathcal{C}^n(A,M):=\text{Hom}(A^{\otimes n},M)$
and that $\mathcal{C}^n(A,M):=M\otimes A^{\otimes n}$.

Consider a symmetric algebra $A$.
By definition, it means that the algebra $A$
is equipped with an isomorphism
$\lineaire:A\buildrel{\cong}\over\rightarrow A^{\vee}$ of $A$-bimodules.
By functoriality,
$\mathcal{C}^*(A,\lineaire):
\mathcal{C}^*(A,A)\buildrel{\cong}\over\rightarrow
\mathcal{C}^*(A,A^{\vee})$
is an isomorphism of cosimplicial modules.
The adjunction map~(\ref{adjunction map})
$Ad:\mathcal{C}^*(A,A^{\vee})
\buildrel{\cong}\over\rightarrow
\mathcal{C}_*(A,A)^{\vee}$ is an isomorphism of
cosimplicial modules (Compare with~\cite[1.5.5]{LodayJ.:cych}).
Let $t_n:\mathcal{C}_n(A,A)\rightarrow\mathcal{C}_n(A,A)$
be the cyclic operator defined by~\ref{cyclic operator}.
The Hochschild chain complex $\mathcal{C}_*(A,A)$
is a cyclic module~\cite[2.1.0]{LodayJ.:cych}.
So $\mathcal{C}_*(A,A)^{\vee}$ with $\tau_n:=t_n^{\vee}$ is a cocyclic module.
Therefore by isomorphism, $\mathcal{C}^*(A,A)$ is also a cocyclic module.

Theorem~\ref{Hochschild algebre symmetrique}
claims that this cocyclic structure on $\mathcal{C}^*(A,A)$
defines a structure of Batalin-Vilkovisky on the Gerstenhaber algebra
$HH^*(A,A)$.

\begin{proof}[Proof of Theorem~\ref{Hochschild algebre symmetrique}]
Let $\varphi:A\otimes A\rightarrow\Bbbk$ be a bilinear form on $A$.
It is easy to see that the associated right linear map
$\lineaire:A\rightarrow A^{\vee}$ is a morphism of $A$-bimodules
if and only if $\varphi$ is symmetric and
\begin{equation}\label{associative de la forme bilineaire}
\varphi(a_2,a_0a_1)=\varphi(a_0,a_1a_2),\quad\forall a_0,a_1,a_2\in A.
\end{equation}
Therefore the endomorphism operad of the symmetric algebra $A$,
$\mathcal{E}nd_A$ is cyclic: it is the cyclic endomorphism operad
defined above.
By definition, $\tau_2\mu$ is the unique map $A\otimes A\rightarrow A$
such that $$\varphi(a_0,\tau_2(\mu)(a_1,a_2))=\varphi(a_2,a_0a_1),
\quad\forall a_0,a_1,a_2\in A.$$
Therefore, by~(\ref{associative de la forme bilineaire}),
we have $\tau_2\mu=\mu$.
In the proof of Corollary~\ref{operad d'endomorphisme d'une algebre},
we have seen that $\mathcal{E}nd_A$ is an operad with
multiplication. Therefore, $\mathcal{E}nd_A$ is a cyclic operad with
multiplication and by Theorem~\ref{Theoreme BV algebre}, $HH^*(A,A)$ is a
Batalin-Vilkovisky algebra.
\end{proof}
\section{Cyclic cohomology of Hopf algebras}\label{cyclic cohomology of Hopf algebras}
Let $\mathcal{H}$ be an Hopf algebra with antipode $S$
and unity $\eta:\Bbbk\rightarrow\mathcal{H}$, $\eta(1_\Bbbk)=1$.
Consider a morphism of algebras (called {\it character}) $\Character:\mathcal{H}\rightarrow\Bbbk$.
The {\it twisted antipode} $\widetilde{S}$ is by definition the convolution
product of $\eta\circ\Character$ and $S$ in $\text{Hom}(\mathcal{H},\mathcal{H})$.
Explicitly, for $h\in\mathcal{H}$, $\widetilde{S}(h)=\Character(h^{1})S(h^{2})$,
where $\Delta h=h^{1}\otimes h^{2}$.
\begin{numerotation}
The couple $(\Character,1)$ is a {\it modular pair in involution} for the
Hopf algebra $\mathcal{H}$
if $\widetilde{S}\circ\widetilde{S}=id_{\mathcal{H}}$.
\end{numerotation}
The twisted antipode $\widetilde{S}$ is an algebra antihomomorphism:
$$\widetilde{S}(ab)=\widetilde{S}(b)\widetilde{S}(a),
\quad\forall a,b\in\mathcal{H},
\quad\widetilde{S}(1)=1.$$
It is also a coalgebra twisted antihomomorphism:
$$
\Delta\widetilde{S}(h)=S(h^{2})\otimes\widetilde{S}(h^{1}),
\quad\forall h\in\mathcal{H}.
$$
More generally, we have
\begin{equation}\label{twisted coalgebra}
\forall n\geq 1,\quad
\Delta\widetilde{S}(h)=S(h^{n})\otimes\dots\otimes S(h^{2})\otimes\widetilde{S}(h^{1}).
\end{equation}

Consider the map $\tau_n:\mathcal{H}^{\otimes n}\rightarrow
\mathcal{H}^{\otimes n}$ defined by
\begin{multline*}
\tau_n(h_1\otimes\dots\otimes h_n):=
\left(\Delta^{n-1}\widetilde{S}(h_1)\right)\cdot
(h_2\otimes\cdots\otimes h_n\otimes 1)\\
=\widetilde{S}(h_1)^{1}h_2\otimes\cdots\widetilde{S}(h_1)^{n-1}h_n\otimes
\widetilde{S}(h_1)^{n}.
\end{multline*}
Here $\cdot$ is the product in $\mathcal{H}^{\otimes n}$
and $\Delta^{n-1}\widetilde{S}(h_1)=\widetilde{S}(h_1)^{1}\otimes\cdots
\otimes \widetilde{S}(h_1)^{n}$
(Review the notations introduced in the proof of 
Corollary~\ref{Cobar d'une bialgebre}).

In~\cite{Connes-Mosco:hopfacctit,Connes-Mosco:cyclchas},
Connes and Moscovici has shown that the Cobar construction on
$\mathcal{H}$ equipped with the maps $\tau_n$ is a cocyclic
module if $(\delta,1)$ is a modular pair in involution.
\begin{proof}[Proof of Theorem~\ref{Cobar algebre de Hopf}]
In the proof of Corollary~\ref{Cobar d'une bialgebre}, we have seen that
$\Omega\mathcal{H}$ is an operad with multiplication.
In order to apply Theorem~\ref{Theoreme BV algebre},
we need to see that $\Omega\mathcal{H}$ is
a cyclic operad with multiplication.
Therefore it remains to prove~(\ref{deux cyclic}) and~(\ref{trois cyclic})
and that $\tau_2\mu=\mu$.

{\bf Proof of~(\ref{deux cyclic}).}
Let $(a_1,\cdots,a_m)\in\mathcal{H}^{\otimes m}$
and $(b_1,\cdots,b_n)\in\mathcal{H}^{\otimes n}$.
Since $\widetilde{S}$ is an algebra antihomomorphism
and $\Delta^{m+n-2}$ is an algebra morphism,
$$\Delta^{m+n-2}\widetilde{S}(a^{1}_1 b_1)=
\Delta^{m+n-2}\widetilde{S}(b_1)\cdot
\Delta^{m+n-2}\widetilde{S}(a^{1}_1).$$
\begin{multline*}
\text{So}\quad\tau_{m+n-1}\left[(a_1,\cdots,a_m)\circ_1 (b_1,\cdots,b_n)\right]\\
=\Delta^{m+n-2}\widetilde{S}(b_1)\cdot
\Delta^{m+n-2}\widetilde{S}(a^{1}_1)\cdot
(a^2_1,\cdots,a_1^n,1,\cdots,1)\cdot
(b_2,\cdots,b_n,a_2,\cdots,a_m,1).
\end{multline*}
Since $\Delta$ is coassociative,
$$\left(\widetilde{S}(b_1)^{(1)},\cdots,\widetilde{S}(b_1)^{(n-1)}
,\widetilde{S}(b_1)^{(n)1},\cdots,\widetilde{S}(b_1)^{(n)m}\right)
=\Delta^{m+n-2}\widetilde{S}(b_1).$$
\begin{multline*}
\text{So}\quad\tau_n(b_1,\cdots,b_n)\circ_n\tau_m(a_1,\cdots,a_m)\\
=\Delta^{m+n-2}\widetilde{S}(b_1)\cdot
(1,\cdots,1,\widetilde{S}(a_1)^{1},\cdots,\widetilde{S}(a_1)^{m})\cdot
(b_2,\cdots,b_n,a_2,\cdots,a_m,1).
\end{multline*}
Therefore to prove~(\ref{deux cyclic}), it suffices to prove
that
\begin{equation}\label{par recurrence}
\Delta^{m+n-2}\widetilde{S}(a^{1}_1)\cdot
(a^2_1,\cdots,a_1^n,1,\cdots,1)=
(1,\cdots,1,\widetilde{S}(a_1)^{1},\cdots,\widetilde{S}(a_1)^{m}).
\end{equation}
Since $\widetilde{S}$ is a twisted antihomomorphism of
coalgebras (\ref{twisted coalgebra}),
\begin{multline}\label{formule}
\Delta^{m+n-2}\widetilde{S}(a^{1}_1)\cdot
(a^2_1,\cdots,a_1^n,1,\cdots,1)\\
=(S\left(a_1^{(1)m+n-1}\right)a_1^2,S\left(a_1^{(1)m+n-2}\right)a_1^3,
\cdots,S\left(a_1^{(1)m+1}\right)a_1^n,\\
S\left(a_1^{(1)m}\right),\cdots,
S\left(a_1^{(1)2}\right),\widetilde{S}\left(a_1^{(1)1}\right)).
\end{multline}
We prove~(\ref{par recurrence}) by induction on $n\in\mathbb{N}^{*}$:

Case $n=1$. Since $a^{1}_1=\Delta^{0}a_1=a_1$,
the two terms of~(\ref{par recurrence}) are equal to
$\Delta^{m-1}\widetilde{S}(a_1)$.

Case $n\geq 2$. Suppose that~(\ref{par recurrence}) is true for $n-1$.

$$\Delta^{m+n-2}\widetilde{S}(a^{1}_1)\cdot
(a^2_1,\cdots,a_1^n,1,\cdots,1)$$
using~(\ref{formule}), since $S$ is an antipode
$$=(\varepsilon\left(a_1^{(1)m+n-1}\right)1,
S\left(a_1^{(1)m+n-2}\right)a_1^{2},\cdots,
S\left(a_1^{(1)m+1}\right)a_1^{n-1},$$
$$S\left(a_1^{(1)m}\right),\cdots,S\left(a_1^{(1)2}\right),
\widetilde{S}\left(a_1^{(1)1}\right))$$
$$=(1,
S\left(\varepsilon\left(a_1^{(1)m+n-1}\right)a_1^{(1)m+n-2}\right)a_1^{2},
\cdots,S\left(a_1^{(1)m+1}\right)a_1^{n-1},$$
$$S\left(a_1^{(1)m}\right),\cdots,S\left(a_1^{(1)2}\right),
\widetilde{S}\left(a_1^{(1)1}\right))$$
since $\varepsilon$ is a counit
$$=(1,
S\left(a_1^{(1)m+n-2}\right)a_1^{2},
\cdots,S\left(a_1^{(1)m+1}\right)a_1^{n-1},$$
$$S\left(a_1^{(1)m}\right),\cdots,S\left(a_1^{(1)2}\right),
\widetilde{S}\left(a_1^{(1)1}\right))$$
using~(\ref{formule}) with $n$ replaced by $n-1$,
$$=(1,\Delta^{m+n-3}\widetilde{S}(a^{1}_1)
\cdot(a^{2}_1,\cdots,a^{n-1}_1,1,\cdots,1))$$
by induction hypothesis
$$=(1,1,\cdots,1,\widetilde{S}(a_1)^1,\cdots,\widetilde{S}(a_1)^m).$$

{\bf Proof of~(\ref{trois cyclic}).}
Since $\Delta^{n-1}$ is a morphism of algebras,
$$\Delta^{n-1}(\widetilde{S}(a_1)^{i-1}a_i)=
\Delta^{n-1}(\widetilde{S}(a_1)^{i-1})\cdot
\Delta^{n-1}(a_i).$$
\begin{multline*}
\text{So}\quad\tau_m(a_1,\cdots,a_m)\circ_{i-1}(b_1,\cdots,b_n)\\
=(\widetilde{S}(a_1)^{1}a_2,\cdots,\widetilde{S}(a_1)^{i-2}a_{i-1},\\
\Delta^{n-1}(\widetilde{S}(a_1)^{i-1})\cdot
\Delta^{n-1}(a_i)\cdot (b_1,\cdots,b_n),\\
\widetilde{S}(a_1)^{i}a_{i+1},\cdots,\widetilde{S}(a_1)^{m-1}a_m,
\widetilde{S}(a_1)^{m}).
\end{multline*}
Since $\Delta$ is coassociative (case $n\geq 1$) and counitary
(case $n=0$),
\begin{multline*}
\widetilde{S}(a_1)^{1}\otimes\cdots\otimes\widetilde{S}(a_1)^{i-2}
\otimes\Delta^{n-1}\widetilde{S}(a_1)^{i-1}\otimes
\widetilde{S}(a_1)^{i}\otimes\cdots\otimes\widetilde{S}(a_1)^{m}\\
=\Delta^{m+n-2}\widetilde{S}(a_1).\\
\text{Therefore}\quad\tau_m(a_1,\cdots,a_m)\circ_{i-1}(b_1,\cdots,b_n)\\
=\Delta^{m+n-2}\widetilde{S}(a_1)\cdot
(a_2,\cdots,a_{i-1},
\Delta^{n-1}(a_i)\cdot(b_1,\cdots,b_n),a_{i+1},\cdots,a_m,1)\\
=\tau_{m+n-1}\left((a_1,\cdots,a_m)\circ_i(b_1,\cdots,b_n)\right).
\end{multline*}

The multiplication $\mu$ on the operad $\Omega\mathcal{H}$ is
$1\otimes 1$. Since $\widetilde{S}(1)=1$, it is easy to check that
$\tau_2\mu=\mu$.
\end{proof}

\section{Proof of parts a) and b) of Theorem~\ref{Theoreme BV algebre}}\label{preuve relation BV}
\begin{proof}[Proof of part a) of Theorem~\ref{Theoreme BV algebre}]
 Let $f\in O(n-1)$.
By~(\ref{deux cyclic}) and $\tau_2\mu=\mu$
$$
\tau_n\delta_1f=\tau_n(f\circ_1\mu)=
\tau_2\mu\circ_2\tau_{n-1}f=\delta_0\tau_{n-1}f.
$$
 By~(\ref{trois cyclic}), for $2\leq i\leq n-1$
$$\tau_n\delta_i f=\tau_n(f\circ_i\mu)=
\tau_n f\circ_{i-1}\mu=\delta_{i-1}\tau_{n-1}f.
$$
By~(\ref{deux cyclic}),
$$\tau_n\delta_n f=\tau_n(\mu\circ_1 f)=
\tau_{n-1}f\circ{n-1}\mu=\delta_{n-1}\tau_{n-1}f.
$$
Let $g\in O(n+1)$. By~(\ref{trois cyclic}),
$$\tau_n\sigma_j g=\tau_n(g\circ_{j+1}e)=\tau_n g\circ_j
e=\sigma_{j-1}\tau_n g.$$
Therefore~\cite[6.1.1]{LodayJ.:cych} the cosimplicial module $O$ is in
fact a cocyclic module.
\end{proof}
Denote by $B$ the Connes coboundary map associated to the cocyclic
module $O$. By~\ref{Theoreme Gerstenhaber}, we already know
that
$H(\mathcal{C}^{*}(O))$ is a Gerstenhaber algebra.
Therefore to prove part b) of Theorem~\ref{Theoreme BV algebre},
it suffices to prove that~(\ref{relation BV}) holds in
cohomology.

{\bf Normalization}.
We would like to use the normalized cochain complex instead of the unnormalized
one, since the formula for Connes coboundary map $B$ is simpler in the
normalized cochain complex.
By definition, the normalized cochain complex associated to $O$, denoted
$\overline{\mathcal{C}}^{*}(O)$, is the subcomplex of $\mathcal{C}^{*}(O)$
defined by
$$\overline{\mathcal{C}}^{n}(O):=\left\{f\in\mathcal{C}^{n}(O)\text{ such that }
\sigma_jf=0\text{ for }0\leq j\leq n-1\right\}.$$
It is well known that the inclusion
$\overline{\mathcal{C}}^{*}(O)\buildrel{\simeq}
\over\hookrightarrow\mathcal{C}^{*}(O)$ is a cochain homotopy equivalence.
It is easy to see that if $f\in\overline{\mathcal{C}}^{m}(O)$
and $g\in\overline{\mathcal{C}}^{n}(O)$ then
$f\cup g\in\overline{\mathcal{C}}^{m+n}(O)$ and $f\circ_i g\in\overline{\mathcal{C}}^{m+n-1}(O)$ for $1\leq i\leq m$.
Therefore $\overline{\mathcal{C}}^{*}(O)$ is both a subalgebra and a sub Lie algebra of $\mathcal{C}^{*}(O)$.
And so, it suffices to show that for any cycles $f$ and
$g\in\overline{\mathcal{C}}^{*}(O)$, (\ref{relation BV}) holds modulo coboundaries.

{\bf Reduction}. In this section, we show that in order to
prove~(\ref{relation BV}), it suffices to prove Proposition~\ref{moitie BV}.
The idea behind that reduction is to start by proving the following particular
case of~(\ref{relation BV}) (where the number of terms has been divided by two):
If $f\in H(\mathcal{C}^{*}(O))$ is of even degree
then $B(f\cup f)$ is divisible by $2$ and
$$f\overline{\circ}f=\frac{1}{2}\{f,f\}=\frac{1}{2}B(f\cup f)-(Bf)\cup f.$$
Proposition~\ref{moitie BV} is a slight generalization of this formula.
Lemma~\ref{multiple de deux} implies in particular
that $B(f\cup f)$ is a multiple of $2$ if $f$ is of
even degree.

The bilinear map of degree $-1$
$$Z:\overline{\mathcal{C}}^{m}(O)\otimes
\overline{\mathcal{C}}^{n}(O)\rightarrow\overline{\mathcal{C}}^{m+n-1}(O),
\quad f\otimes g\mapsto Z(f,g)$$
is defined by
\begin{equation}\label{definition de Z}
Z(f,g):=(-1)^{mn}\sum_{j=1}^{m}(-1)^{j(m+n-1)}
\tau^{-j}_{m+n-1}\sigma_{m+n}(g\cup f).
\end{equation}
Here $\sigma_n:O(n)\rightarrow O(n-1)$ is the extradegeneracy operator
defined by $\sigma_n:=\sigma_{n-1}\tau_n$.

In order to prove Lemma~\ref{multiple de deux}, we need the following two equations
\begin{equation}\label{degenerescence du cup}
\sigma_{m+n}(f\cup g)=\tau_m f\circ_m g.
\end{equation}
\begin{proof}
By~(\ref{trois cyclic}), (\ref{deux cyclic}) and since $\tau_2\mu=\mu$,
\begin{multline*}
\tau_{m+n}(f\cup g)=\tau_{m+n}((\mu\circ_1 f)\circ_{m+1}g)
=\tau_{m+1}(\mu\circ_1 f)\circ_m g\\
=(\tau_m f\circ_m\tau_2\mu)\circ_m g
=(\tau_m f\circ_m\mu)\circ_m g.
\end{multline*}
Therefore since $\mu\circ_2 e=id$
$$\sigma_{m+n}(f\cup g)=\left[(\tau_m f\circ_m\mu)\circ_m g\right]
\circ_{m+n}e=\tau_m f\circ_m g.$$
\end{proof}
\begin{equation}\label{tau echange f et g}
\tau_{m+n-1}^{-n}\sigma_{m+n}(f\cup g)=\sigma_{m+n}(g\cup f).
\end{equation}
\begin{proof}
Using~((\ref{deux generalise})) and equation~(\ref{degenerescence du cup}),
$$
\tau_{m+n-1}^{-n}\sigma_{m+n}(f\cup g)=
\tau_{m+n-1}^{-n}(\tau_m f\circ_m g)=\tau^{-n}_n g\circ_n f=
\sigma_{m+n}(g\cup f)
$$
\end{proof}
\begin{proof}[Proof of Lemma~\ref{multiple de deux}]
The operator $N:O(n-1)\rightarrow O(n-1)$ is defined by
$$N:=\sum_{i=0}^{n-1}(-1)^{i(n-1)}\tau^{i}_{n-1}=
\sum_{j=1}^{n}(-1)^{j(n-1)}\tau^{-j}_{n-1}.$$
By definition, Connes normalized cochain coboundary is
$B:=N\sigma_n:O(n)\rightarrow O(n-1)$.
Therefore, using equation~(\ref{tau echange f et g}),
\begin{align*}
B(f\cup g)&=\sum_{j=1}^{n}(-1)^{j(m+n-1)}
\tau^{-j}_{m+n-1}\sigma_{m+n}(f\cup g)\\
&+\sum_{j=n+1}^{m+n}(-1)^{j(m+n-1)}
\tau^{-(j-n)}_{m+n-1}\tau^{-n}_{m+n-1}\sigma_{m+n}(f\cup g)\\
&=\sum_{j=1}^{n}(-1)^{j(m+n-1)}
\tau^{-j}_{m+n-1}\sigma_{m+n}(f\cup g)\\
&+\sum_{j=1}^{m}(-1)^{(j+n)(m+n-1)}
\tau^{-j}_{m+n-1}\sigma_{m+n}(g\cup f)\\
&=(-1)^{mn}Z(g,f)+Z(f,g)
\end{align*}
\end{proof}
\begin{numerotation}\label{definition de H}
Let $f\in\overline{\mathcal{C}}^{m}(O)$
and $g\in\overline{\mathcal{C}}^{n}(O)$.
Define for any $1\leq j\leq p\leq m-1$
$$H_{j,p}(f,g):=(-1)^{jm-j+(n-1)(p+1+m)}
\tau_{m+n-2}^{-j}\sigma_{m+n-1}(f\circ_{p-j+1}g).$$
and consider the bilinear map of degree $-2$,
$$H:\overline{\mathcal{C}}^{m}(O)\otimes
\overline{\mathcal{C}}^{n}(O)\rightarrow\overline{\mathcal{C}}^{m+n-2}(O),$$
$$f\otimes g\mapsto
H(f,g):=\sum_{1\leq j\leq p\leq m-1}H_{j,p}(f,g).$$
\end{numerotation}
\begin{proof}[Proof of part b) of Theorem~\ref{Theoreme BV algebre} assuming Proposition~\ref{moitie BV}]
By applying Proposition~\ref{moitie BV} and Lemma~\ref{multiple de deux},
\begin{align*}
&f\overline{\circ}g+\varepsilon g\overline{\circ}f
+dH(f,g)+\varepsilon dH(g,f)\\
&+H(df,g)+\varepsilon H(dg,f)
+(-1)^{m-1}H(f,dg)+\varepsilon(-1)^{n-1}H(g,df)\\
&=(-1)^{m}Z(f,g)+\varepsilon(-1)^{n}Z(g,f)
-(-1)^{m}(Bf)\cup g-\varepsilon(-1)^{n}(Bg)\cup f\\
&=(-1)^{m}\left(B(f\cup g)-(Bf)\cup g-(-1)^{m(n-1)}(Bg)\cup f\right).
\end{align*}
Here the sign $\varepsilon$ is equal to $-(-1)^{(m-1)(n-1)}=(-1)^{mn+m+n}$. 
Since in cohomology, the cup product is graded commutative,
relation~(\ref{relation BV}) is proved.
\end{proof}

{\bf Proof of Proposition~\ref{moitie BV}}.
Recall that since $f\in O(m)$,
$$df=\mu\circ_2
f+\sum_{i=1}^{m}(-1)^{i}f\circ_i\mu
+(-1)^{m+1}\mu\circ_1 f.$$
It is easy to see that Proposition~\ref{moitie BV} is a consequence of
the following six equations.
\begin{equation}\label{equation un}
(-1)^{m}Z(f,g)-f\overline{\circ}g=
H(\mu\circ_2 f,g)+H\left((-1)^{m+1}\mu\circ_1 f,g\right).
\end{equation}
\begin{equation}\label{equation deux}
\sum_{1\leq j<p\leq m} H_{j,p}((-1)^{p-j}f\circ_{p-j}\mu,g)
=(-1)^{m} H(f,\mu\circ_2 g).
\end{equation}
\begin{equation}\label{equation trois}
\sum_{1\leq j\leq p\leq m-1}
H_{j,p}\left((-1)^{p-j+1}f\circ_{p-j+1}\mu,g\right)
=(-1)^{m}H\left(f,(-1)^{n+1}\mu\circ_1 g\right).
\end{equation}
\begin{equation}\label{equation quatre}
\sum_{1\leq j\leq
  m}H_{j,m}\left((-1)^{m-j+1}f\circ_{m-j+1}\mu,g\right)=
-(-1)^{m}(Bf)\cup g.
\end{equation}
\begin{multline}\label{equation cinq}
\sum_{1\leq j\leq p\leq m}
H_{j,p}(\sum_{\substack{
 1\leq i\leq m,\\
i\neq p-j,i\neq p-j+1 
}}
(-1)^{i}f\circ_i\mu,g)\\
=-\mu\circ_2 H(f,g)-(-1)^{m+n-1}\mu\circ_1 H(f,g)\\
-\sum_{1\leq j\leq p\leq m-1}
\sum_{\substack{1\leq i\leq p-1\text{ or}\\
p+n\leq i\leq m+n-2}}(-1)^{i}H_{j,p}(f,g)\circ_i\mu.
\end{multline}
\begin{equation}\label{equation six}
\sum_{\substack{{1\leq j\leq p\leq m-1}\\ p\leq i\leq p+n-1}}
(-1)^{i}H_{j,p}(f,g)\circ_i\mu=
(-1)^{m}H\left(f,\sum_{i=1}^{n}(-1)^{i}g\circ_i\mu\right).
\end{equation}
\begin{proof}[Proof of~(\ref{equation un})]
By separating the terms $j=p$ and $j<p$,
\begin{align*}
H(\mu\circ_2 f,g)
&=\sum_{1\leq j\leq p\leq m}(-1)^{jm+(n-1)(p+m)}
\tau_{m+n-1}^{-j}\sigma_{m+n}\left((\mu\circ_2
  f)\circ_{p-j+1}g\right)\\
&=(-1)^{m}Z(f,g)\\
&+\sum_{1\leq j<p\leq m}(-1)^{jm+(n-1)(p+m)}
\tau_{m+n-1}^{-j}\sigma_{m+n}(id\cup f\circ_{p-j}g).
\end{align*}
On the other hand, since (\ref{tau echange f et g})
$\tau_{m+n-1}^{-1}\sigma_{m+n}(f\circ_{p-j+1}g\cup id)=
\sigma_{m+n}(id\cup f\circ_{p-j+1}g),$
\begin{align*}
&(-1)^{m+1}H(\mu\circ_1 f,g)\\
&=\sum_{1\leq j\leq p\leq m}
(-1)^{m+1+jm+(n-1)(p+m)}
\tau_{m+n-1}^{-(j-1)}\sigma_{m+n}(id\cup f\circ_{p-j+1}g).
\end{align*}
Therefore, since (\ref{degenerescence du cup})
$\sigma_{m+n}(id\cup f\circ_{p}g)
=\tau_1 id\circ_1(f\circ_p g)=f\circ_p g$,
by the change of variables $j'=j-1$,
\begin{multline*}
(-1)^{m+1}H(\mu\circ_1 f, g)=-f\overline{\circ} g\\
-\sum_{1\leq j'<p\leq m} (-1)^{j'm+(n-1)(p+m)}
\tau^{-j'}_{m+n-1}\sigma_{m+n}(id\cup f\circ_{p-j'} g).
\end{multline*}
\end{proof}
\begin{proof}[Proof of~(\ref{equation deux})]
By the change of variables $p'=p-1$,
\begin{align*}
&\sum_{1\leq j<p\leq m}H_{j,p}\left((-1)^{p-j}f\circ_{p-j}\mu,g\right)\\
&=\sum_{1\leq j\leq p'\leq m-1}
(-1)^{p'+1-j+jm+(n-1)(p'+1+m)}
\tau^{-j}_{m+n-1}\sigma_{m+n}\left(f\circ_{p'+1-j}(\mu\circ_2 g)\right)\\
&=(-1)^{m}H(f,\mu\circ_2 g).
\end{align*}
\end{proof}
The proof of ~(\ref{equation trois}) is similar.

To prove the last three equations, we will express all the
formulas in terms of composite $\circ_i$ of the elements
$\tau_m^{-j}f$, $g$, $\mu$ and $e$, using again and
again~(\ref{deux generalise}) and~(\ref{trois generalise}).
Therefore, we start by giving a new expression for $H_{j,p}(f,g)$:
\begin{equation}\label{deuxieme formule pour H}
H_{j,p}(f,g)=(-1)^{jm-j+(n-1)(p+1+m)}
\sigma_{j-1}(\tau_{m}^{-j}f\circ_{p+1} g).
\end{equation}
\begin{proof}
We have seen that $O$ is a cocyclic module.
Therefore~\cite[Remark 1.2]{Burghelea-Fiedorowicz:chak},
the following relation between $\tau_n$ and the degeneracy maps $\sigma_i$
holds
$$
\forall\; 0\leq r\leq i\leq n,\quad
\tau_n^{r}\sigma_i=\sigma_{i-r}\tau_{n+1}^{r}.
$$
For the extra degeneracy map $\sigma_{n+1}$, we have
$$
\forall\; 0\leq r\leq n,\quad
\tau_{n}^{r}\sigma_{n+1}=\sigma_{n-r}\tau_{n+1}^{r+1}
$$
or equivalently
\begin{equation}\label{extra degenerescence et tau}
\forall\; 1\leq j\leq {n+1},\quad
\tau_{n}^{-j}\sigma_{n+1}=\sigma_{j-1}\tau_{n+1}^{-j}.
\end{equation}
Therefore using~(\ref{trois generalise}),
$$\tau_{m+n-2}^{-j}\sigma_{m+n-1}(f\circ_{p-j+1}g)
=\sigma_{j-1}\tau^{-j}_{m+n-1}(f\circ_{p-j+1}g)
=\sigma_{j-1}(\tau^{-j}f\circ_{p+1}g). 
$$
\end{proof}
\begin{proof}[Proof of~(\ref{equation quatre})]
By~(\ref{extra degenerescence et tau}),
$$
B(f)=\sum_{j=1}^{m}(-1)^{j(m-1)}\sigma_{j-1}\tau^{-j}_m f.
$$
By~(\ref{deuxieme formule pour H}) and~(\ref{deux generalise}),
\begin{align*}
&\sum_{1\leq j\leq m} H_{j,m}\left((-1)^{m-j+1}f\circ_{m-j+1}\mu,g\right)\\
&=\sum_{j=1}^{m}(-1)^{m+1+jm+j}
\sigma_{j-1}\left[(\tau_2^{-1}\mu\circ_1\tau_m^{-j}f)\circ_{m+1}g\right]\\
&=-(-1)^{m}(Bf)\cup g.
\end{align*}
\end{proof}
\begin{proof}[Proof of~(\ref{equation cinq})]
In all this proof, we put $\varepsilon:=
(-1)^{i+jm+(n-1)(p+m)}$. By~(\ref{deuxieme formule pour H}),
\begin{align*}
&\sum_{1\leq j\leq p\leq m}
H_{j,p}(\sum_{\substack{
 1\leq i\leq m,\\
i\neq p-j,i\neq p-j+1 
}}
(-1)^{i}f\circ_i\mu,g)\\
&=\sum_{\substack{
1\leq j\leq p\leq m,1\leq i\leq m,\\
i+j<p\text{ or }i+j>p+1}
}
\varepsilon\sigma_{j-1}
\left[\tau^{-j}_{m+1}(f\circ_i\mu)\circ_{p+1}g\right].
\end{align*}
Remark that when $m+1\leq i+j$, we can forget the condition
$i+j\neq p$ under theses sums, and that when $m+2\leq i+j$,
we can also forget the condition
$i+j\neq p+1$.

Using respectively (\ref{trois generalise}), (\ref{deux generalise}),
(\ref{deux generalise}) again and
(\ref{deux generalise}) twice, we obtain that 
$$
\tau^{-j}_{m+1}(f\circ_i\mu)=
\begin{cases}
\tau^{-j}_m f\circ_{i+j}\mu &\text{if $i+j\leq m$},\\
\mu\circ_1\tau^{-j}_m f & \text{if $i+j=m+1$},\\ 
\mu\circ_2\tau^{-(j-1)}_m f &\text{if $i+j=m+2$},\\
\tau^{-(j-1)}_m f\circ_{i+j-m-2}\mu &\text{if $m+3\leq i+j$}.
\end{cases}
$$

By the change of variables $i'=i+j$, we have
$\varepsilon=(-1)^{i'-j+jm+(n-1)(p+m)}$,
\begin{align*}
&\sum_{\substack{1\leq j\leq p\leq m,\\1\leq i<p-j}}
\varepsilon\sigma_{j-1}
\left[(\tau^{-j}_m f\circ_{i+j}\mu)\circ_{p+1}g\right]\\
&=\sum_{1\leq j<i'<p\leq m}\varepsilon
\gamma(\tau^{-j}_mf;id,\dots,id,e,id,\dots,id,\mu,id,\dots,id,g,id,\dots,id)\\
&=-\sum_{1\leq j\leq i<p\leq m-1}(-1)^{i}
H_{j,p}(f,g)\circ_i\mu,
\end{align*}
where in the second sum, $e$ is the $j$-th element after the semi-colon,
$\mu$ is the $i'$-th element and $g$ is the $p$-th element.
And we have
\begin{align*}
&\sum_{\substack{1\leq j\leq p\leq m,\\p-j+1< i\leq m-j}}
\varepsilon\sigma_{j-1}
\left[(\tau^{-j}_m f\circ_{i+j}\mu)\circ_{p+1}g\right]\\
&=\sum_{\substack{1\leq j\leq p\leq m,\\
p+1<i'\leq m}}
\varepsilon
\gamma(\tau^{-j}_mf;id,\dots,id,e,id,\dots,id,g,id,\dots,id,\mu,id,\dots,id)\\
&=-\sum_{\substack{
1\leq j\leq p\leq m-1,\\
p+n\leq i\leq m+n-2}}
(-1)^{i} H_{j,p}(f,g)\circ_i\mu,
\end{align*}
where in the second sum, $e$ is the $j$-th element after the semi-colon,
$g$ is the $(p+1)$-th element and $\mu$ is the $i'$-th element.

For $k=1$ or $k=2$,
$$\mu\circ_k H(f,g)=\sum_{1\leq j\leq p\leq m-1}
(-1)^{jm-j+(n-1)(p+1+m)}
\mu\circ_k\left[(\tau^{-j}_m f\circ_{p+1}g)\circ_j e\right].
$$
Since when $1\leq j\leq p\leq m$ and $i+j=m+1$,
we have $1\leq i\leq m$ and the equivalence
$$i+j\neq p+1\Longleftrightarrow p\neq m,$$
$$
\sum_{\substack{1\leq j\leq p\leq m,1\leq i\leq m,\\
i+j=m+1,i+j\neq p+1}}
\varepsilon\sigma_{j-1}\left[(\mu\circ_1\tau_m^{-j}f)\circ_{p+1} g\right]
=-(-1)^{m+n-1}\mu\circ_1 H(f,g).
$$
Since when $1\leq j\leq p\leq m$ and $i+j=m+2$, we have the equivalence
$$
1\leq i\leq m\Longleftrightarrow 2\leq j,
$$
by the change of variables $j'=j-1$ and $p'=p-1$,
$$
\sum_{\substack{1\leq j\leq p\leq m,1\leq i\leq m,\\
i+j=m+2}}
\varepsilon\sigma_{j-1}\left[(\mu\circ_2\tau_m^{-(j-1)}f)\circ_{p+1} g\right]
=-\mu\circ_2 H(f,g).
$$
By the change of variables $j'=j-1$, $i'=i+j-m-2$ and then $p'=p-1$,
we have
$\varepsilon=(-1)^{i'-j'-1+j'm+(n-1)(p+m)}$
and 
\begin{align*}
&\sum_{\substack{1\leq j\leq p\leq m,\\m+3-j\leq i\leq m}}
\varepsilon\sigma_{j-1}
\left[(\tau^{-(j-1)}_m f\circ_{i+j-m-2}\mu)\circ_{p+1}g\right]\\
&=\sum_{1\leq i'<j'<p\leq m}\varepsilon
\gamma(\tau^{-j}_mf;id,\dots,id,\mu,id,\dots,id,e,id,\dots,id,g,id,\dots,id)\\
&=-\sum_{1\leq i'<j'\leq p'\leq m-1}(-1)^{i'}
H_{j',p'}(f,g)\circ_{i'}\mu,
\end{align*}
where in the second sum, $\mu$ is the $i'$-th element after the
semi-colon,
$e$ is the $j'$-th element and $g$ the $p$-th element.
\end{proof}
To prove~(\ref{equation six}), use~(\ref{deuxieme formule pour H})
and the change of variables $i'=i-p+1$.
\section{Proof of part c) of Theorem~\ref{Theoreme BV algebre}}
Using the following Proposition, we see immediately that
part c) of Theorem~\ref{Theoreme BV algebre} follows from parts a) and b).

In this section, all the graded modules are considered as lower graded.
Recall that a {\it mixed complex}
is a graded module $M=\{M_i\}_{i\in\mathbb{Z}}$
equipped with a linear map of degree $-1$ $d:M_i\rightarrow M_{i-1}$ and
a linear map of degree $+1$ $B:M_i\rightarrow M_{i+1}$ such that
$d^{2}=B^{2}=dB+Bd=0$. 

\begin{proposition}\label{complexe mixe et BV-algebre}
Let $(M,d,B)$ be a mixed complex such that
its homology $H_*(M,d)$ equipped with $H_*(B)$ has a Batalin-Vilkovisky
algebra structure. Then its cyclic homology $HC_*(M)$
is a graded Lie algebra of lower degree $+2$.
\end{proposition}
The key point in the proof of this proposition is the following lemma
not explicited stated in~\cite{Chas-Sullivan:stringtop}.
The proof of this lemma is exactly the proof of Theorem 6.1
of~\cite{Chas-Sullivan:stringtop}.
\begin{lem}\label{Chas et Sullivan}
Let $H$ be a Batalin-Vilkovisky algebra and $HC$ be a graded module.
Consider a long exact sequence of the form
$$
\cdots\rightarrow H_n\buildrel{I}\over\rightarrow HC_n
\rightarrow HC_{n-2}\buildrel{\partial}\over\rightarrow H_{n-1}
\rightarrow\cdots
$$If the operator $B:H_i\rightarrow H_{i-1}$ is equal to $\partial\circ
I$ then
$$[a,b]:=(-1)^{\vert a\vert}I\left(\partial a\cup\partial b\right),
\quad\forall a,b\in HC$$
defines a Lie bracket of degree +2 on $HC$.
\end{lem}
\begin{proof}[Proof of Proposition~\ref{complexe mixe et BV-algebre}]
A mixed complex $M$ is a
(differential graded) module over the differential exterior graded algebra
$\Lambda:=(\Lambda\varepsilon_1,0)$~\cite{Kassel:cychcmc}.
Consider the Bar construction of $\Lambda$ with coefficients in $M$,
$B(M;\Lambda;\Bbbk)$.
By~\cite[Proposition 1.4]{Kassel:cychcmc}, the cyclic homology of $M$,
is the homology of $B(M;\Lambda;\Bbbk)$:
$$HC_*(M):=H_*(B(M;\Lambda;\Bbbk))=\text{Tor}^\Lambda_*(M,\Bbbk).$$
Explicitly, $B(M;\Lambda;\Bbbk)$ is the complex defined as follow:

$$B(M;\Lambda;\Bbbk)_n=M_n\oplus M_{n-2}\oplus M_{n-4}\oplus\cdots
\quad\text{and}$$
$$d(m_n,m_{n-2},m_{n-4},\cdots)=(dm_n+Bm_{n-2},dm_{n-2}+Bm_{n-4},\cdots).$$
A mixed complex yields a Connes long exact sequence
$$\cdots\rightarrow H_n(M,d)\buildrel{I}\over\rightarrow HC_n(M)
\buildrel{S}\over\rightarrow HC_{n-2}(M)
\buildrel{\partial}\over\rightarrow H_{n-1}(M,d)\rightarrow\cdots
$$
(Usually $\partial$
is unfortunately denoted by $B$, since it is induced by $B$.)
The connecting homomorphism  $\partial:H_{n-2}(B(M,\Lambda,\Bbbk))\rightarrow H_{n-1}(M,d)$
maps the class of the $n-2$ cycle
$(m_{n-2},m_{n-4},\cdots)$ to
the class of the cycle $Bm_{n-2}$~\cite[Proof of Prop 2.3.6]{Seibt:cychomalg}.
Of course, $I:H_n(M,d)\rightarrow H_{n}(B(M,\Lambda,\Bbbk))$
maps the class of the cycle $m_n$ to $(m_n,0,0,\cdots)$.
Therefore $H_*(B)=\partial\circ I$~\cite[Ex. 9.8.2]{Weibel:inthomalg}.
Finally by Lemma~\ref{Chas et Sullivan},
$HC_*(M)$ is a graded Lie algebra
of degree $+2$.
\end{proof}
\section{Comparison with McClure and Smith}\label{resultat de McClure
  and Smith}
In this section, we compare 
Part b) of Theorem~\ref{Theoreme BV algebre} with a result
announced by McClure and Smith~\cite{McClure-Smith:slidesnorthwestern}:

A cosimplicial module (resp. space) $X^\bullet$
has a {\it cup-cocyclic} structure
if it is a cocyclic module (resp. space) and
if it has a cup product~\cite[Definition
2.1(iii)]{McClure-Smith:deligneconj}
such that
$$
\tau ^{-n-1}_{m+n+1}
(f\cup\delta_0 g)=g\cup\delta_0 f,\quad\forall f\in X^{m},g\in X^{n}. 
$$
\begin{resultat}~\cite[Corollary]{McClure-Smith:slidesnorthwestern} If
the cosimplicial module $X^\bullet$ has a cup-cocyclic structure then
the normalized cochain complex
associated, $\overline{\mathcal{C}}^{*}(X^\bullet)$, has an action by an
operad equivalent to the singular chains on the operad
$\mathcal{D}$ of framed little disks.
\end{resultat}
So $H^{*}(\overline{\mathcal{C}}^{*}(X^\bullet))$ has an action by the operad $H_*(\mathcal{D})$,
i. e.  $H^{*}(\overline{\mathcal{C}}^{*}(X^\bullet))$
is a Batalin-Vilkovisky algebra.
Tedious computations show
that a cosimplicial module (resp. space) has a cup-cocyclic structure
if and only if it is a linear (resp. topological)
cyclic operad with multiplication in our sense.
Therefore their result gives a Deligne's version of our Theorem.

Note that McClure and Smith have announced a topological counterpart
to their result:
\begin{resultat}~\cite[Corollary]{McClure-Smith:slidesnorthwestern} If
the cosimplicial space $X^\bullet$ has a cup-cocyclic structure then
its realisation, $\text{Tot}(X^\bullet)$, has an action by an
operad equivalent to the operad
$\mathcal{D}$ of framed little disks.
\end{resultat}
\bibliography{Bibliographie}
\bibliographystyle{amsplain}
\providecommand{\bysame}{\leavevmode\hbox to3em{\hrulefill}\thinspace}
\providecommand{\MR}{\relax\ifhmode\unskip\space\fi MR }
\providecommand{\MRhref}[2]{%
  \href{http://www.ams.org/mathscinet-getitem?mr=#1}{#2}
}
\providecommand{\href}[2]{#2}

\end{document}